\documentclass[a4paper, twoside, 12pt]{amsart}
\usepackage[margin=3.7cm]{geometry}

\usepackage{times}
\usepackage{graphicx}
\usepackage{framed}

\usepackage[T1]{fontenc}							
\usepackage[latin1]{inputenc}

\usepackage[all]{xy}
\usepackage{pb-diagram, pb-xy}

\usepackage{amssymb}
\usepackage{amsthm}

\newtheoremstyle{myconj}
{0.8em}
{0.4em}
{\itshape}
{}
{\bfseries}
{.}
{ }
{\bfseries \thmname{#1} \thmnumber{#2}\thmnote{(#3)}}%

\newcounter{conjcount}

\newcounter{appAlemcount} 

\newcounter{appBlemcount} 

\theoremstyle{myconj}

\theoremstyle{plain}
\newtheorem*{thm*}{Theorem}
\newtheorem*{lem*}{Lemma}

\newtheorem*{rem*}{Remark}
\newtheorem*{cor*}{Corollary}
\newtheorem*{ex*}{Example}
\newtheorem*{riglem*}{Rigidity lemma}
\newtheorem*{keylem*}{Key lemma}

\newcommand{\bbP}{{\mathbb P}}

\begin{document}
%
%
%
%
%

\title[On subvarieties of Abelian Varieties with  degenerate Gau\ss\ Mapping]{On subvarieties of Abelian Varieties with  degenerate Gau\ss\ Mapping}

\author[Rainer Weissauer]{Rainer Weissauer}

%
%
%
%
%
%

\begin{abstract} Let $X$ be an abelian variety. Then
the {\lq{conormal bundle}\rq}  of an irreducible subvariety $Y$ of $X$ defines an irreducible
Lagrangian subvariety $\Lambda_Y$ of the cotangent bundle $T^*(X)\cong X \times
Lie(X)^*$ of $X$. Our aim is to show that the projection morphism $\gamma: \Lambda_Y \to Lie(X)^*$
is generically finite unless $Y$ is stable under translation by a nontrivial abelian subvariety. 
\end{abstract}

\maketitle

{\bf Introduction}. 
For  a $d$-dimensional irreducible subvariety $Y$ of an abelian variety $X$ of dimension $g$ the Gau\ss\ mapping $\Gamma$ 
attaches to a regular point $y$ of $Y$ the tangent space $T_y(Y)$ of $Y$ at $y$, considered as a $d$-dimensional
linear subspace of the tangent space $T_y(X)$ of $X$ at $y$. If the latter is identified with the Lie algebra $T_0(X)$ by a translation, 
this defines a point of the Gra\ss mann variety $Gr(d,T_0(X))$. For varying $y$, this defines an algebraic morphism on the subset $S\!=Y_{reg}$ of regular points of  $Y$
$$\Gamma\!: S \to Gr(d,T_0(X))\ .$$   Properties of the Gau\ss\
mapping  have been extensively studied by Ochiai, Kawamata, Bogomolov, Ueno, Noguchi, Mori, Abramovich, Ran and others; see e.g. for an overview [D], [D2]. In particular,  Abramovich [A] has shown that $\Gamma$ fails to be a generically finite morphism if and only if the subvariety $Y$ is {\it degenerate} in the sense that there exists an abelian subvariety $A$ of $X$ of dimension $>0$ such that $A+Y=Y$ holds.  

Instead of the above mentioned Gau\ss\ mapping $\Gamma$ involving tangent spaces, 
we consider another type of Gau\ss\ mapping that involves cotangent spaces.
An irreducible subvariety $Y$ of the abelian variety $X$ defines a Lagrangian
subvariety $\Lambda_Y$ of the cotangent bundle $T^*(X)$ of $X$, which for smooth $Y$ is just the conormal bundle of $Y$ in $X$. For an abelian variety the cotangent bundle splits $T^*(X) \cong X \times T_0^*(X)$. Hence the  projection onto the second factor $T_0(X)^*$ (the dual of the Lie algebra) induces an algebraic morphism $$\gamma\!:\! \Lambda_Y \to T_0^*(X)\ .$$ 
The main result of this paper (theorem 1) is the following analog of the result of Abramovich:  If the mapping $\gamma\! :\! \Lambda_Y \to T_0(X)^*$ is not dominant, then $Y$ is degenerate.
For $d\! =\! 1$ and $d\! =\! g-1$ both $\gamma$ and $\Gamma$ can be identified. Otherwise, 
by the obvious identification $Gr(d,T_0(X)) \cong Gr(g-d,T_0(X)^*)$,
the Gau\ss\ mapping $\Gamma$ can also be viewed as a morphism  $S \to Gr(g-d,T_0(X)^*)$, and 
as such may be naturally extended to a morphism $\alpha: \Lambda_Y\vert S \to V$ from the normal bundle $\Lambda_Y\vert S$ over $S$ to  
the tautological vector bundle $V$ of degree $g-d$ over $Gr(g-d,T_0(X)^*)$. 
Since $V$ is a subbundle of the trivial vector bundle $Gr(g-d,T_0(X)^*)\times T_0(X)^*$,
we dispose over a natural projection $\beta\! :\! V \to T_0(X)^*$ such that $\gamma\vert_S = \beta \circ
\alpha$. So clearly, $dim(\Gamma(S)) < d$ implies $\dim(\alpha(\Lambda_Y\vert_S)) < g$.
Hence, by dimension reasons, the mapping $\gamma \! :\!  \Lambda_Y \to T_0(X)^*$ cannot be dominant in case of $dim(\Gamma(S)) < d$.  Thus, our main result is closely related to the geometry of the classical Gau\ss\ mapping and reproves the result of Abramovich for $\Gamma$ (theorem 2). On the other hand, the result of Abramovich for $\Gamma$ does not seem to imply the result for $\gamma$ in a direct way. But, nevertheless, our proof makes significant  use of the methods developed in [A] for the study of the Gauss\ mapping $\Gamma$. 

Our interest in the mapping $\gamma$ comes from the fact, that the generic degrees $d(Y)= deg( \gamma: \Lambda_Y \to
T^*_0(X))$ are strongly related to the 
singularities of $Y$ and the perverse Euler-Poincare characteristic of $Y$. 
Indeed, an irreducible subvariety $Y$ of a complex abelian variety naturally defines a
characteristic subvariety $Ch(Y)$ of the cotangent bundle $T^*(X)$ via the characteristic variety of the $D$-module associated to the middle perverse intersection cohomology sheaf of $Y$ (see below). This
characteristic variety is a finite union of irreducible Lagrangian subvarieties
of $T^*(X)$, hence each of them is of the form $\Lambda_{Y_\nu}$ for
certain irreducible {\lq{characteristic}\rq} subvarieties $Y_\nu$ of $Y$, with multiplicities $m(Y_\nu)$.  As shown in [FK] the 
Euler-Poincare characteristic $$\chi_Y = \sum_i (-1)^i \dim( H^i(Y, IC_Y[d]))$$ of the
intersection cohomology groups $H^i(Y, IC_Y[d]))$ of $Y$ is the sum of the generic degrees $d(Y_\nu)$
of the Lagrangians $\Lambda_{Y_\nu}$
$$  \chi_Y = \sum_\nu m(Y_\nu)d(Y_\nu) \ .$$
Notice, $Y$ itself is always one of the {\lq{characteristic}\rq} subvarieties $Y_\nu$.  All components $Y_\nu \neq Y$ are contained in the singular locus $Y_{sing}=Y \setminus S$. If such $Y_\nu\neq Y$ do occur or not, seems to depend on the structure of the singularities of $Y$ in a rather complicated and not well understood way.

The study of singularities of theta divisors has a long tradition after the work of Andreotti and Mayer [AM]. 
As an application of theorem 1 we show that the theta divisor $Y$ of 
the Jacobian $X$ of a generic curve, although highly singular by Riemann's theorem,
only admits $Y_\nu=Y$ as singular component (theorem 3).  For the theta divisor $Y$ of a principally polarized abelian variety $X$, the structure map $p_Y\! :\! \Lambda_Y \to Y$ is a birational morphism, but there are examples [Kr] where there appear other $Y_\nu$ than $Y$ defining components of $Ch(Y)$. The above relation between the invariants $\chi_Y$ and $d(Y_\nu)$ is also interesting from the
point of view given in [KrW], where it is shown that $\chi_Y$ can be interpreted as the dimension of an irreducible representation $\omega_Y$ of a reductive group $G(Y)$, both canonically attached to $Y$.
In the case of theta divisors, the groups $G(Y)$ define an interesting stratification of the moduli space of principally polarized abelian varieties. In this context, we refer to 
related results on $d(Y)$ from [J], [GM], [SS].   
 
Concerning the definition of $Ch(Y)$, recall that
any $D$-module $K$ on a smooth complex algebraic variety $X$ has an associated characteristic variety $Ch(K)$ that by definition is a subvariety of the cotangent bundle $T^*(X)$ of $X$. If the $D$-module $K$
is holonomic, all  the irreducible components of the characteristic variety $Ch(K)$ 
are subvarieties $\Lambda_{Y_\nu} \subset T^*(X)$ that
are defined by the conormal bundles of certain irreducible subvarieties $Y_\nu$ of $X$. See [KS], [G].
Any middle perverse sheaf $P$ on $X$ defines a holonomic $D$-module $K$
via the Riemann-Hilbert correspondence. Applied for the perverse intersection cohomology sheaf $P= IC_Y[d]$
of $Y$, the results mentioned above for $Ch(Y)$ are a special case of the following formula:
The sum of the degrees 
of the mappings $\gamma\!:\! \Lambda_{Y_\nu} \to T_0^*(X)$ defined by the components of the
characteristic variety $Ch(K)$, with multiplicities, is the Euler-Poincare
characteristic $\chi(P)$  of the perverse sheaf $P$ defined by the solutions of the holonomic $D$-module $K$
via the Riemann-Hilbert correspondence  [FK]. Notice, $\chi(P)=0$ implies that $P$ is translation invariant by an abelian subvariety $A$ of $X$ of dimension $> 0$. For simple abelian varieties $X$
a short proof for this can be found in [KrW]. In general, all known proofs of this result are more complicated. We wonder whether theorem 1 can be helpful to find a short argument as in [KrW] also in the non-simple case.

The paper is organized as follows: In $\S 1$ we formulate the results
and prove theorem 2 and 3. In $\S 2$ we review some results from [A]. Then, in the remaining part of the paper the proof
of theorem 1 is given via  induction on $g=\dim(X)$. 
In $\S 3$ we make some preparations in order to deal with non-simple abelian varieties. Using an argument from [R], we give in $\S 4$ the proof in the case of simple abelian varieties $X$ (lemma 3) formulated in a way that is suitable for the induction argument. We remark, if $X$ is simple and in addition $Y$ is smooth, this is reminiscent of a result of Hartshorne  [H, $\S 4$]. Indeed, then regular 1-forms on $\Lambda_Y$ descend to $Y$ because $\Lambda_Y\to Y$ defines an ample vector bundle unless $Y$ is degenerate. In the remaining part of $\S 4$, by iterated use of lemma 3, we reduce the proof of theorem 1  to the codegenerate case (corollary 3) that finally is treated in $\S 5$.

\goodbreak

\bigskip
{$\S 1$ \bf Notations and preliminary remarks}. Let $X$ be an abelian variety over an algebraically closed field $k$ of dimension $g$.
For a closed irreducible subvariety $Y$ of dimension $d$ in $X$  we let $S$ denote a Zariski dense open subset of its regular locus $Y_{reg}$. 
For the conormal bundle $p_S\!: T_{S}^*(X) \to S$ let $\Lambda_Y$ denote the closure of $T_S^*(X)$ in the cotangent bundle $T^*(X)$ of $X$. The cotangent bundle is trivial for an abelian variety, and  with respect to the corresponding decomposition
$$T^*(X) \cong X \times T^*_0(X)\ $$  the structure  morphism of the bundle $p_X\!: T^*(X)\to X$ is given by the projection on the first factor. 
This defines a structure morphism $p_Y\!: \Lambda_Y \to Y$ by the upper horizontal arrows of the diagram
$$ \xymatrix{ \Lambda_Y  \ar@{^{(}->}[r]^-{i_Y} & T^*(X) \ar[r]^-{p_X} &  X \cr
\Lambda_S \ar@{^{(}->}[u]\ar[rr]^-{p_S}  & & S \ar@{^{(}->}[u] \cr } ,$$
if we take into account that the image
of $p_X\circ i_Y$ is contained in $Y$, since $Y$ is the Zariski closure of $S$ and $\Lambda_Y$ is the Zariski closure of $\Lambda_S$. Notice, 
$\dim(\Lambda_Y)=g$ and $\dim(\Lambda_Y \setminus \Lambda_S) < g$.
For $y\in S$  let $\Lambda_{S,y}$ denote the conormal space $N^*_y(Y)$ of $Y$ in $T^*_y(X)$ at $y$, i.e. the fiber $p_S^{-1}(y)$.  By a translation, $T^*_y(X)$ will always be identifed with its image in $T^*_0(X) =  - y + T^*_y(X)$, so in this sense  $\Lambda_{S,y} \subset T^*_0(X)$.

\medskip
{\it The Gau\ss\ mapping $\gamma$}. 
The 
projection $T^*(X)\cong X \times T^*_0(X) \to T^*_0(X)$ on the second factor restricted to $\Lambda_Y \subset T^*(X)$ induces the Gau\ss\ mapping $$\gamma\!: \Lambda_Y \to T^*_0(X)\ .$$ 
We write $\lambda = (y,\tau)\in \Lambda_Y$, where $y=p_Y(\lambda)\in Y$ and $\tau \in \Lambda_{Y,y}\subset T^*_0(X)$. Then the image under the Gau\ss\ mapping  is the second component $\gamma(\lambda)=\tau\in T^*_0(X)$. For $y\in S$, the  conormal vector $\tau\in T^*_0(X)$, resp. in $T_{Y,y}^*(X) \subset T_y^*(X)$ after translation, annihilates
the tangent space $T_y(Y)$ of $Y$ at $y$.

\medskip
The case $Y\!=\!X$ is exceptional, since in this case the image $\gamma(\Lambda_Y)$ of the Gau\ss\ mapping is contained in $\{0\}$. For $Y\!\neq\! X$, we can remove both the closure of the zero section in $\Lambda_Y$ and the zero section
in $T^*(X)$ to obtain a proper morphism between the associated projective conormal bundles
$$  \bbP\gamma\!: \bbP\Lambda_Y \to \bbP(T^*_0(X)) \ .$$
Obviously for $Y\neq X$, the morphism $\bbP\gamma$ is dominant if and only if $\gamma$ is dominant. This allows to ignore the trivial vector $\tau=0$ in $T^*_0(X)$ 
in subsequent arguments. If the morphism $\gamma$ is not dominant, the image of $\bbP(\Lambda_Y)$ is a closed subvariety of
$\bbP(T^*_0(X))$, hence contained in a hypersurface that is
defined as the zero locus of some nontrivial homogenous polynomial $F$ on $T^*_0(X)$.

\medskip
Since $\dim(\Lambda_Y)=g$, for $Y\neq X$ the following assertions are  equivalent:

\begin{enumerate} 
\item $\dim(\gamma(\Lambda)) < g $
\item The proper morphism $\bbP\gamma: \bbP\Lambda_S \to \bbP(T^*_0(X)) $ is not dominant.
\item The image $\gamma(\Lambda_S)$ of the Gau\ss\ mapping is contained in the zero locus of a nontrivial homogenous polynomial $F$ on $T^*_0(X)$. 
\item $\gamma\! : \Lambda_Y \to T^*_0(X)$ is not dominant.
\item $\gamma\! : \Lambda_Y \to T^*_0(X)$ is not generically finite.
\item For any point $y$ of general position in $S$ there exists a curve  $C$ in $S$ containing $y$, which is contracted by $\gamma$.
\end{enumerate}

Similarly as for irreducible $Y$, one defines the Gau\ss\ mapping $\gamma$
for reducible closed varieties $Y$. In this case the Gau\ss\ mapping is dominant if and only if the Gau\ss\ mapping for one of its irreducible components is dominant.

\medskip An irreducible variety $Y$ in $X$ will be called {\it degenerate}, if
there exists an abelian subvariety $A\subset X$ of positive dimension with the property $A+Y=Y$. (This notion differs from the one used in [R]).
Our main result is 

\medskip
{\bf Theorem 1}. {\it For a closed irreducible subvariety $Y$ of $X$ the following assertions are equivalent
\begin{enumerate}
\item The Gau\ss\ mapping $\gamma\! : \Lambda_Y \to T^*_0(X)$ is not dominant.
\item $Y$ is degenerate. \end{enumerate}
}

\medskip
{\it The Gau\ss \ mapping $\Gamma$}. If we assign to each point $y\!\in\! S$
the tangent space $T_y(Y)$ in $T_y(X)$, this gives rise to the Gau\ss\ mapping $\Gamma$, now with image in the Gra\ss mann variety 
$$ \Gamma\! : S \to Gr(d,T_0(X)) \ .$$
Here again the tangent spac $T_y(Y)$ is considered as a subspace of $T_0(X)$, using translation by $y\in X$. 

\medskip
The following reproves theorem 4 of Abramovich [A]; see also [R], chapter II.

\medskip
{\bf Theorem 2}. {\it For a closed irreducible subvariety $Y$ of $X$ the following holds:
If the Gau\ss\ mapping $\Gamma\! : S \to Gr(d,T_0(X))$ is not generically finite, then the Gau\ss\ mapping $\gamma\! : \Lambda_Y \to T^*_0(X)$ is not dominant and hence
 $Y$ is degenerate.}

\medskip
{\it Proof}. 
If $\Gamma$ is not generically finite, for any point $y$ of $S$ in general position there exists an algebraic curve $C$ containing $y$, which is contracted under $\Gamma$. Then for $y'$ in $C$ we have $T_y(Y)=T_{y'}(Y)$
in $T_0(X)$, hence $\Lambda_{S,y} = N^*_yY = N^*_{y'}Y = \Lambda_{S,y'}$. Since then $\gamma(y,v)=\gamma(y',v)$ holds for all
$v \in \Lambda_{S,y} \subset T^*_0(X)$, the curve $C \times \{ v \} \subset \Lambda_S$ is 
contracted by $\gamma$, for all $v\in N^*_y(Y) \subset T^*_0(X)$. So there exist points in general position contracted by $\gamma$, hence $\gamma\!: \Lambda_S \to T^*_0(X)$ can not dominant by dimension reasons. 
Therefore $Y$ is degenerate by theorem 1. \qed

\medskip
{\bf Theorem 3.} {\it  For the theta divisor $Y$ of the Jacobian  $X$ of a generic regular projective complex curve $C$ the characteristic variety $Ch(Y)$ is irreducible.}

\medskip
{\it Proof}. We can assume that $C$ is not hyperelliptic. Since $dim(Y)=g-1$, we need not distinguish between $\Gamma$ and $\gamma$. So $d(Y)= {2g -2 \choose g-1} $ is the generic degree of $Y$ for the classical Gau\ss\ mapping [GH, p. 360].
Also the perverse Euler-Poincare characteristic $\chi_Y$ of the theta divisor is ${2g -2 \choose g-1}$; see [W, p.273]. Hence $d(Y_\nu)=0$ for every irreducible component  $\Lambda_{Y_\nu}\neq \Lambda_Y$ of $Ch(Y)$, as follows from the formula of [FK]. So all $Y_\nu \neq Y$ are degenerate by theorem 1. If $Ch(Y)$ were not irreducible, therefore
some $Y_\nu\neq Y$ and hence also $Y$ would contain $y + A$ for some $y\in Y$ and some abelian subvariety $A\subset X$ of dimension $>0$. Since the Jacobian of a generic
curve of genus $g$ is an irreducible abelian variety [CG], this would imply $A=X$ which is not possible.   
 \qed

\medskip
{\it Concerning the proof of Theorem 1}. 
To show (a) $\Longrightarrow$ (b) will cover the rest of this paper.
The converse is trivial: Suppose $Y$ is degenerate and $A+Y=Y$ holds. For an abelian variety $A\subset X$ of $\dim(A)>0$,
let $\tilde Y$ denote the image of $Y$ in the quotient $B=X/A$. Notice,
$A+Y\!=\! Y$ implies $T_y(A)\subset T_y(Y)$ and hence  $\Lambda_{Y,y} \subset T^*_0(B)$. Therefore
$\gamma(\Lambda_S) \!=\! \tilde\gamma(\Lambda_{\tilde Y}) \subset T^*_0(B)$, for the corresponding Gau\ss\ mapping $\tilde \gamma$ of $\tilde Y\subset B$. Since
$\dim(\gamma(\Lambda_S)) \! =\! 
\dim(\tilde\gamma(\Lambda_{\tilde S})) \! \leq\! \dim(B)\! < \! \dim(X)$,  the morphism $\gamma\!:\Lambda_Y \to T^*_0(X)$
is not dominant. 

\medskip
We prove the assertion (a) $\Longrightarrow$ (b)  of theorem 1 by induction on the dimension $d$ of $Y$. The case $d=0$ is trivial. So, let us fix some $d>0$. Suppose theorem 1 is
already proven for  irreducible subvarieties $Y'$ of dimension $\dim(Y')\!<\! d$ of an arbitrary abelian variety $X'$. This assumption
will be maintained during the proof almost until the end of the paper. Furthermore,  it is easy to see that for the proof we may assume that $Y$ generates $X$, i.e.
$$ \langle Y \rangle = X \ .$$
Under these assumptions, we then show  that the assertion of theorem 1 also holds for varieties $Y$ of dimension $d$. 
Before we proceed, let us  recall from [A] 
the following 

\bigskip
{$\S 2$ \bf Characterization of degenerate subvarieties}.
For reduced and irreducible subvarieties $Y$ of an abelian variety $X$ define  
$$ Z(Y)=\bigl\{ y\in Y \ \vert \ \exists X'\subset X, X' \mbox{ closed subgroup of } dim(X')>0, y+X' \subset Y\bigr\} \ .$$ 
Then, according to loc. cit. the following holds

\bigskip
{\bf Proposition 1}. { \it If $Y$ is Zariski closed in $X$, then $Z(Y)$ is Zariski closed in $Y$.}

\medskip
{\bf Proposition 2}. { \it If $Y$ is closed in $X$ and $Z(Y)=Y$ holds, $Y$ is degenerate.}

\bigskip
In loc. cit. this is stated in the more general context of semiabelian varieties.

\medskip
{\bf Lemma 1}. {\it Suppose $U$ is a Zariski open dense subvariety of $Y$, and suppose $Y$ is closed in $X$.
Then $Z(U)=U$ implies $Z(Y)=Y$.} 

\medskip
{\it Proof}. Indeed $Z(U) \subset Z(Y)$ by definition, hence  $Y=\overline U =\overline{Z(U)}  =\overline{Z(Y)}$.
Since $Z(Y)$ is Zariski closed in $Y$  by proposition 1, we get $Z=Z(Y)$. \qed

\medskip
{\bf Remark 1}. Keeping lemma 1 and proposition 2 in mind, we may replace $Y$ by some Zariski dense open subset
$U$ of the nonsingular locus $S=Y_{reg}$ of $Y$. For simplicity, we then often 
tacitly write $U=S$ by abuse of notation.

\medskip
{\bf Remark 2}. Suppose $Y$, or a Zariski open dense subset $U$ of $Y$, somehow is written as the union of (not necessarily finitely many) 
subvarieties $F$. Then  $Z(F)=F$ for all these $F$
implies $Z(U)=U$, hence $Z(Y)=Y$.

\bigskip
{$\S 3$ \bf Exact sequences of abelian varieties}. 
1) Let  $X' \subset X$ be a nontrivial abelian subvariety of dimension $<g$. The image  of $Y\subset X$ 
under the quotient mapping $q\!:X\to \tilde X\!=\!X/ X' $ 
will be considered as a closed irreducible subvariety $\tilde Y$ of $\tilde X$,  endowed with the reduced subscheme structure
$$ \xymatrix{ 0 \ar[r] & X' \ar[r]^i & X\ar[r]^q & \tilde X \ar[r] & 0\cr
 &   & Y\ar@{->>}[r]^q\ar@{^{(}->}[u] & \tilde Y \ar@{^{(}->}[u] & \cr} \ .$$
Our  assumption $\langle Y\rangle \!=\! X$ implies $\langle \tilde Y\rangle \!=\! \tilde X$. Hence, $\dim(\tilde Y) >0$ and the fibers $F_{\tilde y}$  of the
morphism $q\!: Y\to \tilde Y$  have dimension $$ \dim(F_{\tilde y}) < d=\dim(Y)\ .$$ 
For $\tilde y\in \tilde Y$, there exists  $y\in Y$ so that
$ F_{\tilde y} \ =\ q^{-1}(\tilde y) \ \subset \ y + X' $.

\medskip
2) For the proof of theorem 1 we may 
replace $X$ by a finite etale covering and $Y$ by its inverse image. This allows to assume
that $X$ splits (non-canonically) into a direct product
$$ X = X' \times \tilde X \ .$$  
Therefore we may tacitly assume that
some splitting of the exact sequence exists, and has been chosen. 
Then  $$ T^*(X)=T^*(X')\times T^*(\tilde X)\quad
\mbox{ and } \quad  F_{\tilde y} =  Y \cap (\tilde y + X') \ .$$

\medskip
3) For (regular) points $y_1,y_2$ in $Y$ with $q(y_1)=q(y_2)$, the fibers $\Lambda_{Y,y_1} = N^*_{y_1}(Y)$ and 
$\Lambda_{Y,y_2} = N^*_{y_2}(Y)$ usually do not coincide.
Let $i: X'\to X$ be the inclusion or more generally any of its translates $i(x') = y+x'$. 
Then we claim

\medskip
{\bf Lemma 2}. {\it There exists a Zariski dense open subset $U$ of the set of regular points of $\tilde S$ of $\tilde Y$, such that for  regular points $y$ of $Y$ in $q^{-1}(U)$ there exists a  canonical exact sequence of vectorspaces 
$$   \xymatrix{0\ \ar[r] & \ \Lambda_{\tilde S,\tilde y}\ \ar[r] & \ \Lambda_{S,y}\   \ar[r]^-{T^*(i)} &  \ \Lambda'_{F_{\tilde y},y}\ \ar[r] & 0 }\ .$$}

\medskip
So for fixed $\tilde y=q(y)$ in $\tilde Y$ with fiber $F_{\tilde y}$ in $X'$, the variation of the conormal spaces $\Lambda_{S,y}$ for $y\in F_{\tilde y}$ is controlled by the variation of the conormal spaces $\Lambda'_{F_{\tilde y},y}$. 
Notice, for a subvariety $Y'$ of a translate of $X'$, we can define $\Lambda_{Y'}$ in $T^*(X)$, and also $\Lambda'_{Y'}$ in $T^*(X')$. The prime index will indicate that the ambient space is a translate of $X'$.

\medskip
{\it Proof of lemma 2}.  Consider  
$0\to T_y(F_{\tilde y}) \to T_y(Y) \to T_{\tilde y}(\tilde Y) \to 0$. This exact sequence of tangent spaces, at a point $y$ where $q$ is  a smooth morphism locally, maps
to $0\to T(X')\to T(X)\to T(\tilde X)\to 0$. Hence, by the snake lemma we get $0\to N'_y(F_{\tilde y}) \to N_y(Y) \to N_{\tilde y}(\tilde Y) \to 0$. 
The exact sequence in our assertion is the dual sequence. This easily shows the claim. \qed

\medskip
Remark: There is also an exact sequence
$ 0 \to \Lambda_{S,y} \to \Lambda_{F_{\tilde y},y} \to T^*_{\tilde y}(\tilde S) \to 0 $.

\medskip
4) Since  $\Lambda_S = \bigcup_{y\in S} \Lambda_{S,y}$, the image $\gamma(\Lambda_Y)\subset T^*_0(X)$ of the Gau\ss\ mapping
 is the Zariski closure of the union
$$\gamma(\Lambda_S) = \bigcup_{y\in S} \gamma(\Lambda_{S,y})\ .$$
Here, of course, $S$ could be replaced by any Zariski dense open subset. We conclude that the image of $\gamma(\Lambda_Y)$ under the linear mapping
$$ T^*(i)\!: T^*_0(X)\to T^*_0(X') $$
 is the Zariski closure of 
$$ T^*(i)\bigl(\bigcup_{y\in S} \gamma(\Lambda_{S,y})\bigr) = \bigcup_{y\in S} T^*(i)\bigl(\gamma(\Lambda_{S,y})\bigr) = \bigcup_{y\in S} \gamma'(\Lambda'_{F_{\tilde y ,y}}) \ $$
where $\gamma'\!: \Lambda'_{F_{\tilde y,y}} \to T^*_0(X')$ denotes
the Gau\ss\ mapping for $X'$, instead of $X$. Indeed, after replacing $S$ by some Zariski open dense subset (also denoted $S$ by abuse of
notation) there exists a commutative diagram
$$ \xymatrix{ \Lambda_{S,y}  \ar@/_9mm/[dd]_{\gamma}\ar[r]^-{T^*(i)}\ar@{^{(}->}[d] &  \Lambda'_{F_{\tilde y},y}  \ar@/^9mm/[dd]^{\gamma'}\ar@{^{(}->}[d]\cr
T^*(X) \ar@{->>}[d] &  T^*(X') \ar@{->>}[d]\cr
T^*_0(X) \ar[r]^-{T^*(i)} &  T_0^*(X') \cr} $$
 Therefore
$$ T^*(i) \bigl(\gamma(\Lambda_Y)\bigr) \ = \ \overline{ \bigcup_{y\in S} \gamma'(\Lambda'_{F_{\tilde y},y} ) } \ .$$

\medskip
5)  Now assume that the Gau\ss\ mapping $$ \gamma: \Lambda_Y \to T^*_0(X) $$ is {\it not dominant}; in addition we assume $Y\neq X$.  
Then using 2),  there exists  
a homogenous polynomial $F\neq 0$ on $T^*_0(X)=T^*_0(X')\oplus T^*_0(\tilde X)$ whose zero locus contains the image
of the Gauss mapping $\gamma$. 
Suppose 
\begin{itemize}
\item $\tau'\neq 0$ in $T^*_0(X')$ is some fixed vector of general position in $T^*_0(X')$.  
\item $\tau'$ is contained in the subvectorspace $                                                 \Lambda'_{F_{\tilde y,y}}$ of $T^*_0(X')$, for $y\in U \subset Y$ in some fixed Zariski dense open subset $U$ of $S$. 
\end{itemize}
\noindent
Then, by  lemma 2, there exists $\tilde\tau$ in $T^*_0(\tilde X)$ such that the vector
$\tau=(\tau',\tilde \tau)$ in $T^*_0(X')\oplus T^*_0(\tilde X)=T^*_0(X)$
is contained in the linear subspace $\Lambda_{S,y}$ of $T^*_0(X)$, and such that furthermore
$$ (\tau',\tilde\tau) + \Lambda_{\tilde S,\tilde y} \ \subset \ \Lambda_{S,y} \ \subset \ \gamma(\Lambda_S) \subset T_0^*(X) \ .$$
So, the polynomial $F$ vanishes on all the vectors $ (\tau',\tilde\tau) + \Lambda_{\tilde S,\tilde y}$. 

\medskip
Notice, $\Lambda_{\tilde S,\tilde y} =  \{0\}  \times W $ holds for some linear subspace $W$ of $V=T^*(\tilde X)$. For $v\in V$ there exists an expansion of $F(\tau',v)$ 
$$ F(\tau',v) = F_{m,\tau'}(v) + F_{m-1,\tau'}(v) + ... + F_{0,\tau'}(v) \ ,$$
where the  $F_{\nu,\tau'}(v)$ are homogenous polynomials of degree $\nu$ on
$V$. We may suppose $\tilde F := F_{m,\tau'} \neq 0$, since otherwise
$F$ would not depend on $v\in V$; and  since $\tau'$ in $T^*(X')$ is of general position by our assumptions, this would  give as a contradiction $F\!=\!0$. 
For any $v\in V$ and any fixed vector $\tilde\tau \in V$ we have (symbolically)
$$ \tilde F(v) \ = \ \lim_{t \to \infty} \ t^{-m} \cdot F(\tau',\tilde\tau + t\cdot v ) \ .$$
Since $F(\tau',\tilde\tau + W)=0$ vanishes, we get: For every $v\in W\subset V$ and any $t\in k^*$ also $t^{-m}\cdot F(\tau',\tilde\tau+ t\cdot v)=0$ vanishes.
We summarize this as follows:
\begin{itemize}
\item $\tilde F\neq 0$ on $V=T^*_0(\tilde X)$
\item $\tilde F= 0$ on $W=\Lambda_{\tilde S,\tilde y} \subset T^*_0(\tilde X)$.
\end{itemize}              
The polynomial $\tilde F$ does not depend on the particular point $y\in U$. It only depends on the fixed decomposition $T^*_0(X)=T^*_0(X')\oplus T^*_0(\tilde X)$ and on the point $\tau'$ in $T^*_0(X')$, where the latter is in general position by assumption. We obtain  

\medskip
{\bf Corollary 1}. {\it Suppose $Y$ is closed and irreducible in $X$ such that the Gau\ss\ mapping $\gamma\! : \Lambda_Y \to T^*_0(X)$ is not dominant. For 
$0\to X'\to X\to \tilde X\to 0$ given with $0< dim(X')<g$, together with  $\tau'$ in $\gamma'(\Lambda'_{F_{\tilde y},y})\subset T^*(i)(\gamma(\Lambda_S))$ of sufficiently general position in $T^*_0(X')$, the Gau\ss\ mapping
$$ \tilde\gamma: \Lambda_{\tilde Y} \to T^*_0(\tilde X) $$ is not dominant for all $\tilde y$ in a Zariski dense open subset of $\tilde S \subset \tilde Y$}.

\medskip
{\it Proof}. We may assume $Y\!\neq\! X$, since otherwise $\tilde Y\!=\!\tilde X$ and the assertion would be trivial. So, for $\tilde y\in \tilde S$, we know that $(\tau',\tilde\tau) + \Lambda_{\tilde S,\tilde y}
\subset \gamma(\Lambda_S)$ for some $\tilde \tau$, as shown in 5) above. For $v\in 
\Lambda_{\tilde S,\tilde y}$ therefore $\tilde F(v)=0$ holds for a fixed
nontrivial polynomial $\tilde F$ on $T^*_0(\tilde X)$, not depending
on $\tilde y$. Hence $\tilde \gamma$ can not be dominant. \qed

\medskip
{\bf Proposition 3}. {\it For $Y$ closed and irreducible in $X$ suppose $$\gamma\! : \Lambda_Y \to T^*_0(X)$$ is not dominant. Furthermore, given $0\to X'\to X\to \tilde X\to 0$ for an abelian subvariety $0\neq X' \subsetneq X$, suppose $$\tilde \gamma\! : \Lambda_{\tilde Y} \to T^*_0(\tilde X)$$ is dominant. Then, for all $\tilde y$ in a Zariski dense open subset of $\tilde S \subset \tilde Y$, none of the Gau\ss\ mappings $$ \gamma'\! : \Lambda'_{F_{\tilde y}} \to T^*_0(X') $$  is dominant (and the same for the irreducible components $Y'$ of these fibers $F_{\tilde y}$). }

\medskip{\it Proof}.  If $\gamma' \!: \Lambda'_{F_{\tilde y}} \to T^*_0(X') $ is dominant for some $\tilde y\in \tilde S \subset \tilde Y$ in general position, there exists a conormal vector $\tau'\neq 0$
of general position in $T^*_0(X')$  such that
$$ \tau' \ \in \ \gamma'(\Lambda'_{F_{\tilde y,y}}) $$
holds for some point $y$ on the fiber $F_{\tilde y}\subset Y$, i.e. $q(y)=\tilde y$. Hence by corollary 1 $\tilde\gamma\! : \Lambda_{\tilde Y} \to T^*_0(\tilde X)$ is not dominant, which contradicts our assumptions. \qed

\medskip
By induction theorem 1 holds for varieties $Y'$ of dimension $<d$. Thus in the situation of proposition 3, we can apply theorem 1 to the irreducible components $Y'$ of the
fibers $F_{\tilde y}$ in $Y$. They have  dimension $\dim(Y')\leq \dim(Y) - \dim(\tilde Y)< d$ for generic $\tilde y$ in $\tilde Y$, since $\tilde Y$ has dimension $>0$ by our assumption $\langle Y\rangle = X$. 
So, in the situation of the induction step of the proof  for theorem 1, after renaming $X'$ by $A$  the last proposition implies

\medskip
{\bf Proposition 4}. {\it Suppose $Y$ is  irreducible of dimension $d=\dim(Y)$, closed in $X$
 and also generates $X$.
Suppose there exists an abelian subvariety $A \subsetneq X$ such that
\begin{itemize}
\item  $ \gamma\!: \Lambda_Y \to T^*_0(X) $ is not dominant.
\item $\tilde \gamma\!: \Lambda_{\tilde Y} \to T^*_0(\tilde X)$ is dominant, for the image $\tilde Y$ of $Y$ in $\tilde X=X/A$.
\end{itemize}
Then $Y$ is degenerate.
}

\medskip
{\it Proof}. The assumptions on the Gau\ss\ mappings imply $A\!\neq\! 0$ and $\tilde Y\! \neq\! \tilde X$. By the other  assumptions $\tilde X\!\neq\! 0$. Hence $\dim(\tilde Y)\!>\! 0$ by $\langle \tilde Y\rangle\! =\! \tilde X$. Then by proposition 3 all the Gau\ss\ mappings
$\gamma'\!: \Lambda'_{Y'} \to T^*_0(X')$ for the irreducible components $Y'$ of the fibers $F_{\tilde y}$ (for $\tilde y$ in general position) are not dominant. By the general induction assumption
underlying the proof of theorem 1, we conclude that for all these $\tilde y$ the components $Y'$ are degenerate. Hence $Z(Y')\!=\!Y'$ holds, and therefore
$Z(U)\!=\!U$ holds for some Zariski dense open subset $U$ of $Y$ by remark 2. Thus
$Y$ is degenerate by proposition 1 and 2. \qed

\medskip
In the situation of the induction step for theorem 1, proposition 3 also implies

\medskip
{\bf Corollary 2}. {\it Suppose $Y$ is irreducible and closed in $X$
of dimension $d$ and generates $X$. Suppose $Y$ is not degenerate and $ \gamma\! : \Lambda_Y \to T^*_0(X) $ is not dominant.
Then for any abelian subvariety $X' \subset X$ with quotient $\tilde X:=X/X'$ and image $\tilde Y$ of $Y$ in $\tilde X$,
the Gau\ss\ mapping $\tilde \gamma\!: \Lambda_{\tilde Y} \to T^*_0(\tilde X)$ is not dominant.
}

\medskip
{\it Proof}. We can assume  
$\dim(X')>0$, so that proposition 4 can be applied. \qed

\medskip

\bigskip
{$\S 4$ \bf Fibers of the Gau\ss\ mapping}.
Suppose $Y\subsetneq X$ is irreducible and closed, and suppose the Gau\ss\ mapping
$$ \gamma: \Lambda_Y \to T^*_0(X) $$
is not dominant. In this section we furthermore assume $\langle Y\rangle\! =\!X$ and $Y\!\neq\! X$. 
In this setting we now use the following argument of [R], or [KrW]: 
First observe that under these assumptions all nonempty fibers of the Gau\ss\ mapping $\gamma$
$$  Z_\tau = \gamma^{-1}(\tau) \subset \Lambda_Y $$
 have dimension
$$ \dim(Z_\tau) \geq 1 \ .$$
This follows from the upper semicontinuity of fiber dimensions, since it holds
for generic points $\tau$ in $\gamma(\Lambda_Y)$. Notice, the image $Y_\tau = p_Y(Z_\tau)$ in $Y$ 
$$ \xymatrix{ Z_\tau \ar@{^{(}->>}[d]_{p_Y}\ar@{^{(}->}[r] & \Lambda_Y \ar[d]^{p_Y}\ar[r]^-\gamma & T^*_0(X) \cr
Y_\tau \ar@{^{(}->}[r] & Y & \cr }$$
has the property that
$ p_Y\!: Z_\tau \to Y_\tau $
is a set theoretic bijection, since over any $y\in Y_\tau$ the points $z\in Z_\tau$ are uniquely determined by the condition $\gamma(z)\!=\!\tau$.  Indeed, set theoretically,
$$ Z_\tau = Y_\tau \times \{ \tau \} \subset X \times T^*_0(X) = T^*(X) \ .$$
Now assume $\tau\neq 0$. Then $\tau$ defines a nontrivial linear form
$$ \tau\!: T_0(X) \to k $$
whose kernel contains all tangent vectors in $T_y(Y)$, considered as vectors in $T_0(X)$ via a translation by $y\in X$.
In particular, all tangent vectors in $T_y(Y_\tau)$ at regular points $y$ of $Y_\tau$ are contained in the kernel of the linear form $\tau$.

\medskip
So, for given $\tau\!\neq\! 0$, let us fix some point $y\!=\!y_\tau$ in $Y_\tau$. 
Then the translate $Y_\tau - y_\tau$ contains zero, and the abelian subvariety $X'$ of $X$
generated by $Y_\tau - y_\tau$ is a nontrivial abelian subvariety $X'\subsetneq X$. Indeed, $\dim(Y_\tau)\geq 1$ implies $X'\!\neq\! 0$ and $\tau(T_0(X'))=0$ implies $X'\neq X$. Furthermore by construction of $X'$ $$ Z_\tau \subset y_\tau + X' \ .$$
In this situation, a priori, the abelian variety $X'=X'(\tau,y)$ may depend on the choice of $\tau\in \gamma(\Lambda)$ and also on the choice of $y=y_\tau$ in $Y_\tau$, so that 
$$ \xymatrix{ Y_\tau  \ar@{^{(}->}[r]\ar@{^{(}->}[d]  & X  \ar[r]\ar@{^{(}->}[d] & T^*_0(X) \ni \tau \cr
F_{\tilde y_\tau} = Y\cap \bigl(y_\tau + X'(\tau,y_\tau)\bigr) \ar@{^{(}->}[r] & X \cr} $$
and
$$ Y = \bigcup_{\tau \in \gamma(\Lambda_Y)} Y_\tau \ .$$

\medskip
{\bf Rigidity Property}. {\it For all $\tau$ in a Zariski open dense subset of $\gamma(\Lambda_Y)$ and all $y_\tau$ in a Zariski open dense subset of $Y_\tau$, the abelian variety $X'=X'(\tau,y_\tau)$ does not depend on the choice of $\tau$ and $y_\tau$. 
}

\medskip
{\it Proof}. There exist only countably many abelian subvarieties $X'$ in $X$, and $X'=X'(\tau,y_\tau)$ depends algebraically on $\tau$ and $y_\tau$. Replace $k$ by an uncountable extension field. \qed

\medskip
By the rigidity property we can assume that $X'$ is a fixed nontrivial proper abelian subvariety attached to $Y\subset X$, so that for all $y$ in a Zariski dense open subset $U\subset Y$ 
$$ F_{\tilde y} = Y \cap (\tilde y + X') = Y \cap (y_\tau + X'(\tau,y_\tau)) $$ contains $Y_\tau$ and hence is of positive dimension
$ dim(F_{\tilde y}) \geq  1 ;$ 
and the image $\tilde Y$ of $Y$ in $\tilde X = X/X'$ is irreducible of dimension $$ \dim(\tilde Y) < \dim(Y) \ .$$
To summarize, this shows

\medskip
{\bf Lemma 3}. {\it For irreducible  closed $Y\!\neq\! X$  with $\langle Y\rangle\! =\!X$ and non-dominant Gau\ss\ mapping $\gamma\!: \Lambda_Y \to T^*_0(X)$  there exits an abelian subvariety $0\neq X' \subsetneq X$  such that
$\dim(\tilde Y) < \dim(Y)$ holds for the image $\tilde Y$ of $Y$ in $\tilde X=X/X'$, with the fibers of the Gau\ss\ mapping $\gamma$ contained in translates of $X'$.}  

\medskip
Indeed $\langle Y\rangle\! =\!X$ can
be assumed without restriction of generality.

\medskip
In the situation of lemma 3, we now assume that $Y$ is not degenerate with a non-dominant Gau\ss\ mapping $\gamma$,
and let us also assume that theorem 1 holds for varieties of dimension $\!<\! \dim(Y)$. Then proposition 4 can be applied;  it shows that the induced Gau\ss\ mapping
$$ \tilde\gamma\!: \tilde Y \to T^*_0(\tilde X) $$
is not dominant. Then, $\langle \tilde Y \rangle \!=\! \tilde X$ is inherited from
$\langle Y \rangle\! =\! X$. So suppose
$$ \tilde Y \neq \tilde X \ .$$
Then, if $\tilde Y\!\neq\! \tilde X$, we can apply 
once again lemma 3, now for the pair $(\tilde Y,\tilde X)$, to construct an exact sequence 
$0\to \tilde X' \to \tilde X \to \tilde{\tilde X} \to 0$ such that
\begin{itemize}
\item $\langle \tilde{\tilde Y} \rangle = \tilde{\tilde X}$
\item ${\tilde{\tilde \gamma}}\!: \tilde{\tilde Y}  \to T^*_0(\tilde{\tilde X} )$
is not dominant.
\end{itemize}
Obviously, this construction can be iterated and terminates after finitely many steps since $$ \cdots < \dim(\tilde{\tilde Y}) < \dim(\tilde Y) < \dim(Y)\ , $$ thus provides an abelian subvariety $AÊ\subsetneq X$ containing $X'$, so that the image of $Y$ in $B=X/A$ is equal to $B$. 

\medskip
A closed irreducible variety $Y$ in $X$ will be called {\it codegenerate} (with respect to $A$), if there exists
an abelian subvariety $A\!\neq\! X$ in $X$ such that the image of $Y$ in $B=X/A$ is equal to $B$. 

\medskip
Using this notion 
in the context of the induction step for the proof of theorem 1, we have now shown in this situation  that for $Y$ closed irreducible of dimension $\dim(Y)=d$  the following corollary holds.

\medskip
{\bf Corollary 3}. {\it  Suppose the Gau\ss\ mapping $$\gamma\!: \Lambda_Y \to T^*_0(X)$$ is not dominant, then either 
\begin{enumerate}
\item $Y$ is degenerate, or 
\item $Y$ is codegenerate with respect to an abelian subvariety $A$, so that the fibers of the Gau\ss\ mapping $\gamma$ are contained in translates of $A$.
\end{enumerate}} 

\medskip 

\medskip
{$\S 5$ \bf The proof of theorem 1 in the codegenerate case}.
To complete the proof of the induction step (hence the proof of theorem 1)
for $Y$ of dimension $d$ with non-dominant Gau\ss\ mapping,
it remains to consider the codegenerate case $\tilde Y=B$
of the last corollary 3.  Of course,  without restriction of generality we can assume  in addition  that $Y$ is not degenerate. That $Y$ is not degenerate implies (by the induction assumption of theorem 1 and proposition 1 and 2):
For a Zariski dense open subset $U$ of $\tilde Y=B$, the Gau\ss\ mappings of the fibers  $F_{\tilde y}, \tilde y\in U$ of the projection $q\!: Y \to \tilde Y=B$  are nondegenerate. Indeed, if this non-degeneracy holds for a single fiber $F_{\tilde y}$ where $\tilde y$ is supposed to be in general position, it holds for all fibers $F_{\tilde y}$ with $\tilde y$ in a Zariski dense open subset $U$ of $\tilde Y$ by a specialization argument. 

\medskip
Assuming these conditions all together, we  
claim:  $Y$ {\it is degenerate}.
This gives a contradiction which implies the induction step for the proof of theorem 1.

\medskip
Recall that in the last section we found
an exact sequence 
$$ \xymatrix{ 0 \ar[r] & A \ar[r]^i & X \ar[r]^q & B \ar[r] & 0} $$
with $B=\tilde Y$ such that the fibers $Z_\tau $ of the Gau\ss\ mapping $\gamma: \Lambda_Y \to T^*_0(X)$ map bijectively to
varieties $Y_\tau \subset Y$ that are contained 
in the fibers $$ F_{\tilde y} = Y \cap (\tilde y + A) \subset Y $$
of the projection $q\!:Y \to \tilde Y=B$. Here, without restriction of generality, we  assume 
that $B$ splits, so that $B$ can be considered as a subvariety of $X$ complementary to $A$.  Since $Y=\bigcup_{\tilde y \in B} F_{\tilde y}$, the variety $Y$ is degenerate if all the $F_{\tilde y}$ are degenerate for all $\tilde y$ in some Zariski open dense subset of $\tilde Y$ (using proposition 1 and 2). Therefore, if $Y$ were not degenerate, by the induction assumption of theorem 1 we conclude that the Gau\ss\ mapping $$ \gamma_A\! : \Lambda_{F_{\tilde y}}^A \to T^*_0(A) $$
is dominant for all points $\tilde y$ of general position in $\tilde Y$. 

\medskip
Fix some $\tau' \!\neq\! 0$  in $T^*_0(A)$ in general position;
fix some $\tilde y$, now with $\tilde y\in U$, so that $\gamma_A\!:  \Lambda_{F_{\tilde y}} \to T^*_0(A)$ is dominant. Since $\gamma_A$
is dominant and since $\tau$ has general position in $T^*_0(A)$, there exists
$$ y \in F_{\tilde y} \subset \tilde y + A \mbox{ so that } 
 \tau'  \mbox{ is contained in } N^*_y(F_{\tilde y}) = \Lambda_{F_{\tilde y},y} \ .$$
Since $\tilde S$ is Zariski dense in $\tilde Y=B$, we get in $T^*(B) = B\times T^*_0(B)$
$$ \Lambda_{\tilde S, \tilde y} = \Lambda_{B,\tilde y} = \{ \tilde y \} \times \{ 0\} \ .$$
By lemma 2 we therefore obtain

\medskip
{\bf Lemma 4}. {\it For $\tilde y$ in a suitably chosen Zariski open dense subset  $U$  of $\tilde Y=B$, our assumptions imply that $\tau'$ (chosen in general position) is in $ \Lambda'_{F_{\tilde y},y}$ such that
$$ \Lambda_{S,y} \cong  \Lambda'_{F_{\tilde y},y} \ .$$}

\medskip
In other words: $\tau'\in T^*_0(A)$  can be uniquely lifted to a point in $\Lambda_{S,y}$, once
the corresponding base point $y\in F_{\tilde y}$ over $\tilde y\in U$ has been specified.

\medskip  For the  conormal bundle in $T^*(A)$ we now write $\Lambda^A_{F_{\tilde y}}$ instead of $\Lambda'_{F_{\tilde y}}$. Notice that $\gamma_A\!: \Lambda^A_{F_{\tilde y}} \to T^*_0(A)$ is dominant and hence generically finite, for all $\tilde y$  a Zariski open subset $\tilde U$ of $\tilde Y$.  We therefore obtain

\medskip
{\bf Lemma 5}. {\it For fixed $\tau'\!\neq\! 0$ with general position in $T^*_0(A)$ and fixed $\tilde y$ with general position in $ B$, there exist only finitely many points $y\in Y$ mapping to $\tilde y$ (i.e. $y \in F_{\tilde y}$) for which $\tau'$ is contained in the conormal bundle $\Lambda_{F_{\tilde y},y}$ of $F_{\tilde y}$ at $y$. 
}

\medskip
Combining lemma  4 and 5 we get

\medskip
{\bf Lemma 6}. {\it We have the following diagram 
$$ \xymatrix{  \mbox{dense open subset of } \gamma(\Lambda_S) \ar@{=}[dd] \ar@{^{(}->}[rr] & & T^*_0(X) \ar[dd]^{T^*(i)}\cr
& & \cr
\bigcup_{\tilde y \in U} \gamma'(\Lambda^A_{F_{\tilde y}}) \ar[rr] & & T^*_0(A)  } $$
and the image under the lower horizontal morphism of each $\gamma'(\Lambda^A_{F_{\tilde y}})$ for $ \tilde y\in U$  is Zariski dense in $T^*_0(A)$.
Furthermore, for $\tau'$ of sufficiently general position in $T^*_0(A)$
the points $\gamma(\lambda) \in T^*_0(X)$ in $\gamma(\Lambda_S)$, which map under $T^*(i): T^*_0(X)\to T^*_0(A)$ to the point $\tau'$,
correspond to points $\lambda$ in $\Lambda_S$
$$ \lambda = (y,\tau) \in \Lambda_S \subset X \times T^*_0(X) $$
for which the base point 
$$ \tilde y = q(y) \in B =X/A $$
of $y$ can be arbitrarily prescribed within a dense open subset of $B$. Once this base point $\tilde y \in B$ is fixed, there exists at least one, but at most finitely many choices for the point $y\in X$ such that
$$ \tilde y = q(y) \mbox{ and } \lambda = (y,\tau)=(y,\tau',\tilde\tau) \in \Lambda_S $$
holds for some $\tilde\tau$ in $T^*_0(B)$.  
}

\medskip
Hence there exists a Zariski dense open subset $V\subset \Lambda_S$, such that on $V$ the mapping $\varphi=(q \circ p_S, T^*(i)\circ \gamma)$ 
defines a generically finite and hence dominant morphism $$ \xymatrix{ \Lambda_S \supset V  \ar[r]^-\varphi & \ B \times T^*_0(A) \cr }  $$
$$ \lambda =(y,\tau) \mapsto (\tilde y, \tau') \ .$$

\medskip
For a point $\tau\in \gamma(\Lambda_S)$ in general position
consider the fiber $Z_\tau \subset \Lambda_S$. Now $V\cap Z_\tau$ is a Zariski dense open subset of $Z_\tau$ by the choice of $\tau$.
Furthermore, for all points $\lambda$
in $V\cap Z_\tau$ by definition $\lambda=(y,\tau) $ holds; and the fixed $\tau$ maps to a fixed
$\tau'$ of general position in $T^*(A)$. In other words: The second component  $T^*(i)\circ \gamma$ of the morphism $\varphi$
is constant on $Z_\tau$. Since $\varphi$ is finite on $V$, therefore the first component $$   q\circ p_S: V \cap Z_\tau \to B $$ of the morphism $\varphi$ 
is a finite morphism on $V\cap Z_\tau$, and
$$  V \cap Z_\tau \ni (y,\tau) \mapsto \tilde y=q(y)\ .$$ 
Since $X' \subset A$, using the notations of the argument preceding corollary 3, we already know from the beginning of this section that
$$ Z_\tau \subset y_\tau + X' \subset y_\tau + A  = \tilde y_\tau + A \ .$$
Here without restriction of generality we assumed $X\!=\!A\times B$, and $\tilde y_\tau \in B$ is the image of $y_\tau$ under the projection $q\!:X\to B\!=\!X/A$. Recall that $y_\tau$ was some fixed chosen point in $Y_\tau$ and only depends on $\tau$.
Hence the inclusion $ Z_\tau \subset \tilde y_\tau + A $
implies 
$$ q\circ p_S(V\cap Z_\tau) = \tilde y_\tau \ .$$
Since $q\circ p_S$ is finite on $V\cap Z_\tau$,
the intersection
$V\cap Z_\tau$ therefore contains only finitely many point. But this contradicts
$\dim(V\cap Z_\tau)\!=\!\dim(Z_\tau)\! \geq\! 1$
and finishes the proof of the induction step for
theorem 1. \qed 

\medskip

\medskip

\goodbreak

\bigskip\noindent

\bigskip\noindent

\centerline{\scshape Mathematisches Institut}
\centerline{\scshape Ruprecht-Karls-Universit\"at Heidelberg}
\centerline{{\it E-mail address:} weissauer@mathi.uni-heidelberg.de}


\medskip

\begin{thebibliography}{XXXX}
 \bibitem[A]{A} Abramovich D., {\em Subvarieties of semiabelian varieties}.  Comp. Math. 90, 37--52 (1994)
\bibitem[AM]{AM} Andreotti A., Mayer A.L.  {\em On the period relations for abelian integrals on algebraic curves}. Ann. Sc. Norm. Super. Pisa Cl. Scr. 21, 257 -- 261 (1967)
 \bibitem[CG]{CG} Ciliberto C., van der Geer G., {\em Subvarieties of the moduli space of curves parametrizing Jacobians wit non-trivial endomorphisms}. Am. J. of Math. vol. 114. no. 3, 551 -- 570 (1992)
 \bibitem[D]{D} Debarre O., {\em Complex tori and abelian varieties}. SMF/AMS Texts and monographs, vol. 11 (1999)
 \bibitem[D2]{D2} Debarre O., {\em Fulton-Hansen and Barth-Lefschetz theorems for subvarieties of abelian varieties}.
 J. reine angew. Math. 467, 187 -- 197 (1995)
 \bibitem[FK]{FK} Franeki J., Kapranov M., {\em The Gauss map and a noncompact Riemann-Roch formula for constructible sheaves on semiabelian varieties}. Duke Math. J. 104, no. 1,
171--180 (2000) 
\bibitem[G]{G} Ginsburg V., {\em Characteristic varieties and vanishing cycles}. Invent. math. 84, 327 -- 402 (1986)
\bibitem[GH]{GH} Griffiths Ph., Harris J., {\em Principles of Algebraic Geometry}. John Wiley \& Sons, (1978)  
\bibitem[GM]{GM} Grushevsky M., Salvati Manni R., {\em Singularities of the theta divisor at points of order 2}. Int. Math. Res. Not. 15 (2007)
\bibitem[H]{H} Hartshorne R., {\em Ample vector bundles on curves}. Nagoya Math. J. 43, 73 -- 89 (1971)
\bibitem[J]{J} de Jong R., {\em Theta functions on the theta divisor}. Rocky Mountain J.
Math. 40, 155--176 (2010)
\bibitem[KS]{KS} Kashiwara M., Shapira P., {\em Micro-hyperbolic systems}. Acta Math. 142, 1 -- 55 (1979)
\bibitem[Kr]{Kr} Kr\"amer T., {\em Cubic threefolds, Fano surfaces and the monodromy of the Gauss map}. arXiv: 1501.00226 
\bibitem[KrW]{KrW} Kr\"amer T.,Weissauer R., {\em Vanishing theorems for constructible sheaves on abelian varieties}. J. Algebraic Geometry 24, 531 -- 568 (2015)
\bibitem[R]{R} Ran Z., {\em On Subvarieties of Abelian Varieties}. Invent. math. 62,
459 -- 479 (1981)
\bibitem[SS]{SS} Smyth B., Sommese A.J., {\em On the degree of the Gauss mapping of a submanifold of an Abelian variety}.
Comment. Math. Helvetici 59, 341--346 (1984)
\bibitem[W]{W} Weissauer R., {\em Tannakian categories attached to abelian varieties}. In: Modular Forms on Schiermonnikoog, edited by Edixhofen B., van der Geer G., Moonen B., Cambridge university press, 267 --274 (2008).

\end{thebibliography}
\end{document}